\documentclass[12pt]{amsart}

\begin{document}
\def\refname{References}
\newtheorem{thm}{Theorem}[subsection]
\newtheorem{lem}[thm]{Lemma}
\newtheorem{cor}[thm]{Corollary}
\newtheorem{claim}[thm]{Claim}

\title{Parametrizing the abstract Ellentuck theorem}
\author{Jos\'e G. Mijares}
\address{Escuela de Matem\'aticas. Facultad de Ciencias, Universidad Central de Venezuela. jmijares@euler.ciens.ucv.ve}
\begin{abstract}
We give a parametrization with perfect subsets of $2^{\infty}$ of
the abstract Ellentuck theorem (see \cite{CaS}, \cite{Tod} or
\cite{Tod2}). The main tool for achieving this goal is a sort of
parametrization of an abstract version of the Nash-Williams
theorem. As corollaries, we obtain some known classical results
like the parametrized version of the Galvin-Prikry theorem due to
Miller and Todorcevic \cite{Mil}, and the parametrized version of
Ellentuck's theorem due to Pawlikowski \cite{Paw}. Also, we obtain
parametized vesions of nonclassical results such as Milliken's
theorem \cite{Milk}, and we prove that the family of perfectly
Ramsey subsets of $2^{\infty}\times FIN_{k}^{[\infty]}$ is closed
under the Souslin operation.
\end{abstract}

\maketitle

\section{Introduction.}

In \cite{Ell}, Ellentuck considers the space
$\mathbb{N}^{[\infty]}$, of all the infinite sets of natural
numbers, with the \textit{exponential} (or \textit{Ellentuck's})
\textit{topology} to obtain his famous topological generalization
of Ramsey's theorem \cite{Ram}. The basic open sets of this
topology are the neighborhoods $[a, \alpha] = \{\beta\in
\mathbb{N}^{[\infty]} : a\subset\beta\subseteq\alpha\}$, where $a$
(resp. $\alpha$) is a finite (resp. infinite) subset of
$\mathbb{N}$. We recall that a set $\mathcal{X}\subseteq
\mathbb{N}^{[\infty]}$ is \textit{Ramsey} (\textit{completely
Ramsey} in \cite{Ell}) if for every Ellentuck neighborhood $[a,
\alpha]$ there exists an infinite $\beta\subseteq\alpha$ such that
either $[a, \beta]\subseteq \mathcal{X}$ or $[a, \beta]\cap
\mathcal{X} = \emptyset$. In \cite{Ell}, Ellentuck proved that a
set is Ramsey if and only if it has the Baire property relative to
Ellentuck's topology, proving in this way that the family of
Ramsey sets is closed under the Souslin operation and hence
obtaining a simpler topological proof of Silver's theorem
\cite{Sil} that classical analytic subsets of
$\mathbb{N}^{[\infty]}$ are Ramsey, and generalizing the
Galvin-Prikry theorem \cite{GP} that classical Borel subsets of
$\mathbb{N}^{[\infty]}$ are Ramsey.

Theorems analog to that of Ellentuck have also been proven on
other spaces and different contexts, provided an analog
"exponential" topology is given in each case (see for instance
\cite{DiT}, \cite{Far} and \cite{Milk}). These are the so called
{\em Ellentuck type theorems}. Spaces where an Ellectuck type
theorem can be proven are called (topological) \textit{Ramsey
spaces}. The main tool for proving such theorems (including
Ellentuck's) is deviced for each particular case as a sort of
"combinatorial forcing", inspired by the works of Galvin-Prikry
\cite{GP} and Nash-Williams \cite{Nas}. Recently, based on a
previous work of Carlson and Simpson (see \cite{CaS}), S.
Todorcevic has isolated abstract combinatorial features which are
conditions of sufficiency for topological spaces provided of a
suitable "exponential" topology, to guarantee that Ellentuck type
theorems can be obtained in such spaces. These combinatorial
features enable to device a sort of "abstract Galvin-Prikry
machinery" which constitutes the main means for obtaining an
abstract Ellectuck theorem. It turns out that spaces satisfying
these conditions are Ramsey spaces. We present a summary of
Todorcevic's recent presentation of the abstract topological
Ramsey theory in section 2.

In this work we will show that the conditions proposed in
\cite{Tod}, and summarized in the next section, are also
sufficient to obtain a parametrized version of the abstract
Ellentuck theorem. In this way, we obtain as corollaries some
known classical results like the parametrized version of the
Galvin-Prikry theorem due to Miller and Todorcevic \cite{Mil}, and
the parametrized version of the Ellentuck theorem due to
Pawlikowski \cite{Paw} which makes use of the notion of abstract
Baire property of Morgan \cite{Mor}. Also, we obtain parametized
vesions of nonclassical results such as Milliken's theorem
\cite{Milk}, and we prove that the family of perfectly Ramsey
subsets of $2^{\infty}\times FIN_{k}^{[\infty]}$ is closed under
the Souslin operation, and hence that analytic subsets of
$2^{\infty}\times FIN_{k}^{[\infty]}$ are perfectly Ramsey. Here
$2^{\infty}$ is the space of all the infinite sequences of 0's and
1's, with the product topology regarding $2 = \{0, 1\}$ as a
disdrete space. The space $FIN_{k}^{[\infty]}$, of infinite block
basic sequences of functions $p : \mathbb{N}\rightarrow \{0,
1,\cdots, k\}$ of finite domain and with $p(n) = k$ for some $n$,
will be described in section \ref{cases}.

\section{Abstract Topological Ramsey Theory. Summary of main facts.}

All the definitions and results throughout this section were taken
from \cite{Tod2} and are expected to appear in \cite{Tod}. A
previous presentation of the following notions can also be found
in \cite{CaS}.

We will consider triplets of the form $(\mathcal{R}, \leq, (p_{n})_{n\in\mathbb{N}} )$,
where $\mathcal{R}$ is a set, $\leq$ is a quasi order on $\mathcal{R}$, $\mathbb{N}$ is
 the set of natural numbers, and for every $n\in\mathbb{N}$, $p_{n}: \mathcal{R}\rightarrow \mathcal{P}_{n}$
 is a function with range $\mathcal{P}_{n}$. For each $A\in \mathcal{R}$, we say that $p_{n}(A)$ is
 \textit{the} $n$th \textit{approximation of} $A$. In order to capture the combinatiorial structure
 required to ensure the provability of an Ellentuck type theorem, we will impose some assumptions on
 $(\mathcal{R}, \leq, (p_{n})_{n})$. The first three of them are the following:

\begin{itemize}
\item[{(A1)}]For any $A\in \mathcal{R}$, $p_{0}(A) = \emptyset$.
\item[{(A2)}]For any $A,B\in \mathcal{R}$, if $A\neq B$ then $(\exists n)p_{n}(A)\neq p_{n}(B)$.
\item[{(A3)}]If $p_{n}(A) = p_{m}(B)$ then $n = m$ and $(\forall i<n)p_{i}(A) = p_{i}(B)$.
\end{itemize}
These three assumptions allow us to identify each $A\in \mathcal{R}$ with the sequence $(p_{n}(A))_{n}$ of its approximations. In this way, if we consider the space $\mathcal{P}=\bigcup_{n} \mathcal{P}_{n}$ with the discrete topology, we can identify  $\mathcal{R}$ with a subspace of the (metric) space $\mathcal{P}^{[\infty]}$ (with the product topology) of all the sequences of elements of $\mathcal{P}$. Via this identification, we will regard $\mathcal{R}$ as a subspace of $\mathcal{P}^{[\infty]}$, and we will say that $\mathcal{R}$ is {\it metrically closed} if it is a closed subspace of $\mathcal{P}^{[\infty]}$.

Also, for $a\in \mathcal{P}$ we define the \textit{lenth} of $a$, $|a|$, as the unique $n$ such that $a = p_{n}(A)$ for some $A\in \mathcal{R}$.

\vspace{.25 cm}

We also consider on $\mathcal{R}$ the \textit{Ellentuck type neighborhoods} $$[a,A] = \{B\in \mathcal{R} : (\exists n)(a = p_{n}(B))\ \ \mbox{and} \ \ B\leq A\}$$
where $a\in\mathcal{P}$ and $A\in \mathcal{R}$. If $[a,A]\neq\emptyset$ we will say that $a$ is \textit{compatible} with $A$ (or $A$ is compatible with $a$). Let $\mathcal{P}[A] = \{a\in \mathcal{P} : a\ \ \mbox{is compatible with}\ \ A\}$.

\vspace{.25 cm}

We write $[n,A]$ for $[p_{n}(A),A]$, and $Exp(\mathcal{R})$ for
the family of all the neighborhoods $[n,A]$. This family generates
the natural "exponential" topology on $\mathcal{R}$ which is finer
than the product topology.

\vspace{.25 cm}

A sequence $([n_k,A_k])_k$ of elements of $Exp(\mathcal{R})$ is
called a \textit{fusion sequence} if it is infinite and if
\begin{itemize}
\item[{(i)}]$(n_k)_k\subseteq\mathbb{N}$ is nondecreasing and
$\lim_{k\rightarrow\infty} n_k = \infty$,
\item[{(ii)}]$A_{k+1}\in
[n_k,A_k]$ for all $k$.
\end{itemize}

Since $\mathcal{R}$ is closed, building fusion sequences
constitute a very useful procedure for defining desired elements
of $\mathcal{R}$: for every fusion sequence $([n_k,A_k])_k$ we
have that $\bigcap_k [n_k,A_k]\neq\emptyset$ (and a singleton).
The \textit{limit} of the fusion sequence is the unique element
$A_{\infty}$ of $\bigcap_k [n_k,A_k]$. Note that
$p_{n_k}(A_{\infty}) = p_{n_k}(A_k)$, for all $k$.

\vspace{.25 cm}

With this notation, we can define an analog notion for subsets of
$\mathcal{R}$, to that of \textit{Ramseyness} for subsets of
$\mathbb{N}^{[\infty]}$:

{\bf Definition 1.0.1.} A set $\mathcal{X}\subseteq \mathcal{R}$
is \textit{Ramsey} if for every neighborhood $[a,A]\neq\emptyset$
there exists  $B\in [a,A]$ such that $[a,B]\subseteq \mathcal{X}$
or $[a,B]\cap \mathcal{X} = \emptyset$. A set
$\mathcal{X}\subseteq \mathcal{R}$ is \textit{Ramsey null} if for
every neighborhood $[a,A]$ there exists  $B\in [a,A]$ such that
$[a,B]\cap \mathcal{X} = \emptyset$.

{\bf Definition 1.0.2.} We say that $(\mathcal{R}, \leq,
(p_{n})_{n})$ is a \textit{Ramsey space} if subsets of
$\mathcal{R}$ with the Baire property are Ramsey and meager
subsets of $\mathcal{R}$ are Ramsey null.

\vspace{.25 cm}

In \cite{Tod} it is shown that A1, A2 and A3, together with the
following three assumptions are conditions of suficiency for a
triplet $(\mathcal{R}, \leq, (p_{n})_{n})$, with $\mathcal{R}$
metrically closed, to be a Ramsey space (see also \cite{CaS}):

\vspace{.25 cm}

(A4)(\textit{Finitization}) There is a quasi order $\leq_{fin}$ on
$\mathcal{P}$ such that:
    \begin{itemize}
    \item[{(i)}]$A\leq B$ iff
    $\forall n\exists m \ \ p_{n}(A)\leq_{fin} p_{m}(B)$.
    \item[{(ii)}]$\{b\in \mathcal{P} : b\leq_{fin} a\}$ is finite, for every
    $a\in \mathcal{P}$.
    \end{itemize}

\vspace{.25 cm}

Given $a$ and $A$, we define the \textit{depth of} $a$ \textit{in}
$A$, $depth_{A}(a)$, as the minimal $n$ such that $a\leq_{fin}
p_{n}(A)$.

\begin{itemize}

\item[{(A5)}](\textit{Amalgamation}) Given compatible $a$ and $A$ with $depth_{A}(a)=n$, the following holds:
\begin{itemize}
\item[{(i)}] $\forall B\in [n,A]\ \ ([a,B]\neq\emptyset)$.
\item[{(ii)}] $\forall B\in [a,A]\ \ \exists A'\in [n,A]\ \ ([a,A']\subseteq [a,B])$.
\end{itemize}
\item[{(A6)}](\textit{Pigeon Hole Principle}) Given compatible $a$
and $A$ with $depth_{A}(a) = n$, for each partition
$\phi:\mathcal{P}_{|a|+1}\rightarrow \{0,1\}$ there is $B\in
[n,A]$ such that $\phi$ is constant in $p_{|a|+1}[a,B]$.

\end{itemize}

\vspace{1 cm}

{\bf Abstract Ellentuck Theorem:}

\vspace{.25 cm}

\begin{thm}\label{AbsEll}
[Carlson] Any $(\mathcal{R}, \leq, \mathcal{P}_{n}, (p_{n})_{n})$
with $\mathcal{R}$ metrically closed and satisfying A1-A6 is a
Ramsey space.
\end{thm}

\qed

\vspace{.25 cm}

For instance, take  $\mathcal{R} = \mathbb{N}^{[\infty]}$, the set
of infinite subsets of $\mathbb{N}$, $\leq \ \ = \ \ \subseteq$
and $p_{n}(A)$ = the first $n$ elements of $A$, for each $A\in
\mathbb{N}^{[\infty]}$ . So, the set of approximations is
$\mathcal{P} = \mathbb{N}^{[<\infty]}$, the set of finite subsets
of $\mathbb{N}$. The family of neighborhoods $[a,A]$, with
$a\in\mathbb{N}^{[<\infty]}$ and $A\in \mathbb{N}^{[\infty]}$, is
the family of Ellentuck neighborghoods defined in the
introduction. Define  $\leq_{fin}$ as $a\leq_{fin} b$ iff
($a=b=\emptyset$ or $a\subseteq b$ and $max(a)=max(b)$), for
$a,b\in\mathbb{N}^{[<\infty]}$. With these definitions, A1-A6 are
easily verified. In this case A6 reduces to a natural variation of
the classical pigeon hole principle for finite partitions of an
infinite set of natural numbers. Note also that
$\mathbb{N}^{[\infty]}$ is easily identified with a closed
subspace of $\mathcal{P}^{[\infty]}$, namely, the set of all the
sequences $(x_{n})_{n}$ of finite sets such that $x_{n} =
x_{n+1}\setminus \{max(x_{n+1})\}$, for each $n\in\mathbb{N}$.
Then $(\mathbb{N}^{[\infty]}, \subseteq, (p_{n})_{n}))$ is a
Ramsey space in virtue of the abstract Ellentuck theorem. Hence,
Ellentuck's theorem is obtained as corollary:

\vspace{.25 cm}

\begin{cor}\label{Ell}
[Ellentuck] Given $\mathcal{X}\subseteq\mathbb{N}^{[\infty]}$, the
following hold:
\begin{itemize}
\item[{(a)}]$\mathcal{X}$ is Ramsey iff $\mathcal{X}$ has the Baire Property, relative to Ellentuck's topology.
\item[{(b)}]$\mathcal{X}$ is Ramsey null iff $\mathcal{X}$ is meager, relative to Ellentuck's topology.
\end{itemize}
\end{cor}

\section{Parametrizing with perfect sets.}

In this section we will present our main result. We recall that
$2^{\infty}$ denotes the space of infinite sequences of 0's and
1's, with the product topology regarding $2 = \{0, 1\}$ as a
discrete space. Also, $2^{<\infty}$ denotes the set of finite
sequences of 0's and 1's. Let us consider some features of the
perfect subsets of $2^{\infty}$, following \cite{Paw}:

Some notation is needed. For $x = (x_{n})_n \in 2^{\infty}$, $x|_{
k}$ denotes the finite sequence $(x_{0}, x_{1}, \cdots, x_{k-1})$.
For $u\in 2^{<\infty}$, let $[u] = \{x\in 2^{\infty}: (\exists k)
(u = x|_{k})\}$ and let $|u|$ be the lenth of $u$. Given a perfect
set $Q\subseteq 2^{\infty}$, let $T_{Q}$ be its asociated perfect
tree.  Also, for $u, v = (v_{0}, v_{1}, \cdots, v_{|v|-1})\in
2^{<\infty}$ we write $u\sqsubseteq v$ to mean $(\exists k\leq
|v|) (u = (v_{0}, v_{1}, \cdots, v_{k-1}))$. For each $u\in
2^{<\infty}$, let $Q(u) = Q\cap [u(Q)]$, where $u(Q)\in T_{Q}$ is
define inductively, as follows: $\emptyset(Q) = \emptyset$.
Suppose $u(Q)$ defined. Find $\sigma\in T_{Q}$ such that $\sigma$
is the $\sqsubseteq$-extension of $u(Q)$ where the first
ramification occurs. Then, set $(u*i)(Q)$ = $\sigma*i$, $i = 0,1$.
Here "$*$" denotes "concatenation". Note that for each $n$, $Q$ =
$\bigcup\{Q(u) : u\in 2^{n}\}$.

Given $n\in \mathbb{N}$, and perfect sets $S$ and $Q$ we say that
$S\subseteq_{n} Q$ if $S(u)\subseteq Q(u)$, for every $u\in
2^{n}$. The relation "$\subseteq_{n}$" is a partial order. If for
every $u\in 2^{n}$ we have chosen a $S_{u}\subseteq Q(u)$, then $S
= \bigcup_{u}S_{u}$ is perfect and we have $S(u) = S_{u}$ and
$S\subseteq_{n} Q$. As pointed out in \cite{Paw}, the most
important feature of this partial order is the \emph{property of
fusion}: if $Q_{n+1}\subseteq_{n+1} Q_{n}$, $n\in \mathbb{N}$,
then the fusion $Q = \bigcap_{n} Q_{n}$ is a perfect set and
$Q\subseteq_{n} Q_{n}$, for each $n$.

\vspace{.5 cm}

Our goal in this section is proving a parametrized version of the
abstract Ellentuck theorem of the previous section. This
constitutes the main result contained in this work and is stated
as thereom \ref{thm1} below.

We introduce now the abstract version of the notion of "perfectly
Ramsey" (see \cite{Paw}). We will use the same name: given a
triplet $(\mathcal{R}, \leq, (p_{n})_{n})$ as defined at the
previous section, we say that a set $\mathcal{X}\subseteq
2^{\infty}\times \mathcal{R}$ is \textit{perfectly Ramsey} if for
every perfect set $Q\subseteq 2^{\infty}$ and every neighborhood
$[a,A]\neq\emptyset$ there exist a perfect set $S\subseteq Q$ and
$B\in [a,A]$ such that $S\times [a,B]\subseteq \mathcal{X}$ or
$S\times [a,B]\cap \mathcal{X} = \emptyset$. A set
$\mathcal{X}\subseteq 2^{\infty}\times \mathcal{R}$ is
\textit{perfectly Ramsey null} if for every perfect set
$Q\subseteq 2^{\infty}$ and every neighborhood
$[a,A]\neq\emptyset$ there exist a perfect set $S\subseteq Q$ and
$B\in [a,A]$ such that $S\times [a,B]\cap \mathcal{X} =
\emptyset$.

Also, we will need to generalize the notion of abstract Baire
property (see \cite{Mor}) to this context:

Let $\mathbb{P}$ be the family of perfect subsets of $2^{\infty}$.
We will say that a set $\mathcal{X}\subseteq 2^{\infty}\times
\mathcal{R}$ has the $\mathbb{P}\times
Exp(\mathcal{R})$-\textit{Baire property} if for every perfect set
$Q\subseteq 2^{\infty}$ and every neighborhood $[a,A]$ there exist
a perfect set $S\subseteq Q$ and a neighborhood $[b,B]\subseteq
[a,A]$ such that $S\times [b,B]\subseteq \mathcal{X}$ or $S\times
[b,B]\cap \mathcal{X} = \emptyset$. A set $\mathcal{X}\subseteq
2^{\infty}\times \mathcal{R}$ is $\mathbb{P}\times
Exp(\mathcal{R})$-\textit{meager} if for every perfect set
$Q\subseteq 2^{\infty}$ and every neighborhood $[a,A]$ there exist
a perfect set $S\subseteq Q$ and a neighborhood $[b,B]\subseteq
[a,A]$ such that $S\times [b,B]\cap\mathcal{X} = \emptyset$.

\vspace{.25 cm}

Now we are ready to state our main result:

\begin{thm}\label{thm1}
Given $(\mathcal{R}, \leq,(p_{n})_{n})$,  with $\mathcal{R}$
metrically closed and satisfying A1-A6, the following are true:

\begin{itemize}
\item[{(a)}]$\mathcal{X}\subseteq 2^{\infty}\times \mathcal{R}$ is perfectly Ramsey iff $\mathcal{X}$ has the $\mathbb{P}\times Exp(\mathcal{R})$-Baire Property.
\item[{(b)}]$\mathcal{X}\subseteq 2^{\infty}\times \mathcal{R}$ is perfectly Ramsey null iff $\mathcal{X}$ is $\mathbb{P}\times Exp(\mathcal{R})$-meager.
\end{itemize}
\end{thm}

\vspace{1 cm}

The main tool for proving theorem \ref{thm1} is the following fact
which is a sort of parametrization of an abstract version of
Nash-Williams' theorem:

\begin{thm}\label{thm2}

Given $(\mathcal{R}, \leq,(p_{n})_{n})$,  with $\mathcal{R}$
metrically closed and satisfying A1-A6, the following is true:

For every $\mathcal{F}\subseteq 2^{<\infty}\times\mathcal{P}$,
perfect $P\subseteq 2^{\infty}$ and $A\in\mathcal{R}$ there exist
a perfect $S\subseteq P$ and $D\leq A$ such that one of the
following holds:

\begin{itemize}
\item[{(a)}]for every $x\in S$ and every $C\leq D$ there exist
integers $k$ and $m > 0$ such that $(x|_{k},p_{m}(C))\in
\mathcal{F}$. \item[{(b)}]$(T_{S}\times
\mathcal{P}[D])\cap\mathcal{F} = \emptyset$.
\end{itemize}

\end{thm}

\vspace{1 cm}

Theorem \ref{thm2} is inspired on a parametrized version of the
semiselective Nash-Williams theorem proved by I. Farah (see
theorems 2.2 and 2.3 of \cite{Far}). Before proving it, we need to
start up our parametrized combinatorial machinery, based on the
tecniques used in \cite{Far}, \cite{Paw}[13] and \cite{Tod2}.

\vspace{.5 cm}

{\bf Combinatorial Forcing 1.} Fix $\mathcal{F}\subseteq
2^{<\infty}\times\mathcal{P}$. For a perfect $Q\subseteq
2^{\infty}$, $A\in \mathcal{R}$ and a pair $(u,a)\in
2^{<\infty}\times\mathcal{P}$, we say that $(Q,A)$
\textit{accepts} $(u,a)$ if for every $x\in Q(u)$ and for every
$B\in [a,A]$ there exist integers $k$ and $m$ such that
$(x|_{k},p_{m}(B))\in \mathcal{F}$. We say that $(Q,A)$
\textit{rejects} $(u,a)$ if for every perfect $S\subseteq Q(u)$
and every $B\leq A$, compatible with $a$, $(S,B)$ does not accepts
$(u,a)$. Also, we say that $(Q,A)$ \textit{decides} $(u,a)$ if it
accepts or rejects it.

\vspace{.25 cm}

{\bf Combinatorial Forcing 2.} Fix $\mathcal{X}\subseteq
2^{\infty}\times\mathcal{R}$. For a perfect $Q\subseteq
2^{\infty}$, $A\in \mathcal{R}$ and a pair $(u,a)\in
2^{<\infty}\times\mathcal{P}$, we say that $(Q,A)$
\textit{accepts} $(u,a)$ if $Q(u)\times [a,A]\subseteq
\mathcal{X}$. We say that $(Q,A)$ \textit{rejects} $(u,a)$ if for
every perfect $S\subseteq Q(u)$ and every $B\leq A$, compatible
with $a$, $(S,B)$ does not accepts $(u,a)$. And as before, we say
that $(Q,A)$ \textit{decides} $(u,a)$ if it accepts or rejects it.

\vspace{.5 cm}

{\bf Note:} Lemmas \ref{lem1}, \ref{lem2}, and \ref{lem3} below
hold {\em for both} combinatorial forcings defined above.
\begin{lem}\label{lem1}
The following are true:
\begin{itemize}
\item[{(a)}]If $(Q,A)$ accepts (rejects) $(u,a)$ then $(S,B)$ also
accepts (rejects) $(u,a)$,  for every perfect $S\subseteq Q(u)$
and every $B\leq A$ compatible with $a$. \item[{(b)}]If $(Q,A)$
accepts (rejects) $(u,a)$ then $(Q,B)$ also accepts (rejects)
$(u,a)$, for every $B\leq A$ compatible with $a$. \item[{(c)}]For
all $(u,a)$ and $(Q,A)$ such that $A$ is compatible with $a$,
there exist a perfect $S\subseteq Q$ and $B\leq A$, compatible
with $a$, such that $(S,B)$ decides $(u,a)$. \item[{(d)}]If
$(Q,A)$ accepts $(u,a)$ then $(Q,A)$ accepts $(u,b)$ for every
$b\in p_{|a|+1}[a,A]$. \item[{(e)}]If $(Q,A)$ rejects $(u,a)$ then
there exist $B\in [depth_A(a),A]$ such that $(Q,A)$ does not
accept $(u,b)$ for every $b\in p_{|a|+1}[a,B]$.
\item[{(f)}]$(Q,A)$ accepts (rejects) $(u,a)$ iff $(Q,A)$ accepts
(rejects) $(v,a)$, for every $v\in 2^{<\infty}$ with $u\sqsubseteq
v$.
\end{itemize}
\end{lem}

\begin{proof}

\vspace{.25 cm}

(a), (b), (c), (d) and (f) follow from the definitions. Now to
proof (e), take $(u,a)$ with $|a| = m$ and suppose $(Q,A)$ rejects
it. Define $\phi : \mathcal{P}_{m+1}\rightarrow 2$ such that
$\phi(b) = 1$ iff $(Q,A)$ accepts $(u,b)$. Let $n = depth_{A}(a)$.
By A6, there exists $B\in [n,A]$ such that $\phi$ is constant on
$p_{m+1}[a,B]$.

If $\phi$ takes value 1 on $p_{m+1}[a,B]$ then $(Q,B)$ accepts $(u,a)$. So, in virtue of part (b), $\phi$ must take value 0 on $p_{m+1}[a,B]$ since $(Q,A)$ rejects $(u,a)$. Then $B$ is as required.

\end{proof}

\vspace{.25 cm}

\begin{lem}\label{lem2}
For every perfect $P\subseteq 2^{\infty}$ and $A\in \mathcal{R}$
there exist a perfect $Q\subseteq P$ and $B\leq A$ such that
$(Q,B)$ decides $(u,a)$, for every $(u,a)\in 2^{<\infty}\times
\mathcal{P}[B]$ with $depth_{B}(a)\leq |u|$.
\end{lem}
\begin{proof}

\vspace{.25 cm}

Let $<>$ be the empty sequence of 0's and 1's, and recall from
section 2 that $\emptyset\in\mathcal{P}$. Using lemma
\ref{lem1}(c), find a perfect $Q_{0}\subseteq P$ and $B_{0}\leq A$
such that $(Q_{0},B_{0})$ decides  $(<>,\emptyset)$.

Suppose we have defined $Q_{n}$ and $B_{n}$ such that
$(Q_{n},B_{n})$decides every $(u,a)\in 2^{n}\times
\mathcal{P}[B_{n}]$ with $depth_{B_{n}}(a) = n$. Let $u_{0}$,
$u_{1}$,... , $u_{2^{n+1}-1}$ be a list of the elements of
$2^{n+1}$; and let $b_{0}$, $b_{1}$,... , $b_{r}$ be a list of the
$b\in \mathcal{P}[B_{n}]$ such that $depth_{B_{n}}(b) = n+1$.

Using lemma \ref{lem1}(c), find a perfect $Q_{n}^{0,0}\subseteq
Q_{n}(u_{0})$ and $B_{n}^{0,0}\leq B_{n}$ compatible with $b_{0}$
such that $(Q_{n}^{0,0},B_{n}^{0,0})$ decides $(u_{0},b_{0})$.  We
can suppose  $B_{n}^{0,0}\in [b_{0},B_{n}]$. Hence, we can assume
$B_{n}^{0,0}\in [n+1,B_{n}]$ in virtue of A5(ii) and lemma
\ref{lem1}(b).

In a similar way, for every $(i,j)\in \{0, 1,\cdots ,
2^{n+1}-1\}\times \{0, 1,\cdots , r\}$, we can find $Q_{n}^{i,j}$
and $B_{n}^{i,j}$ with: $Q_{n}^{i,j+1}\subseteq
Q_{n}^{i,j}(u_{i})$, $B_{n}^{i,j+1}\in [b_{j+1},B_{n}^{i,j}]$,
$Q_{n}^{i+1,0}\subseteq Q_{n}(u_{i+1})$,  $B_{n}^{i+1,0}\in
[b_{0},B_{n}^{i,r}]$; and such that $(Q_{n}^{i,j},B_{n}^{i,j})$
decides $(u_{i},b_{j})$. (Notice that this construction is
possible in virtue of A5(i). Again, we can  assume $B_{n}^{i,j}\in
[n+1,B_{n}]$ in virtue of A5(ii) and lemma \ref{lem1}(b)).

Let $Q_{n+1} = \bigcup_{i = 0}^{2^{n+1}-1}Q_{n}^{i,r}$ and
$B_{n+1} = B_{n}^{2^{n+1}-1,r}$. Then, $(Q_{n+1},B_{n+1})$ decides
$(u,b)$, for every $(u,b)\in 2^{n+1}\times \mathcal{P}[B_{n+1}]$
with $depth_{B_{n+1}}(b) = n+1$: for such $(u,b)$, there exist
$(i,j)\in \{0, 1,\cdots ,  2^{n+1}-1\}\times \{0, 1,\cdots , r\}$
such that $u = u_i$ and $b = b_j$. Then
$(Q_{n}^{i,j},B_{n}^{i,j})$ decides $(u,b)$. Notice that
$$Q_{n+1}(u_i) = Q_n^{i,r}\subseteq Q_n^{i,r-1}(u_i)\subseteq
\dots \subseteq Q_n^{i,j}(u_i)$$ and $$B_{n+1} =
B_{n}^{2^{n+1}-1,r} \leq B_{n}^{i,j}.$$ Hence $(Q_{n+1},B_{n+1})$
decides $(u_i,b_j)$.

\medskip

Besides, $Q_{n+1}\subseteq_{n+1} Q_{n}$ and $B_{n+1}\in
[n+1,B_{n}]$

\medskip

Now let $Q = \bigcap_{n} Q_{n}$ and take $B\in
\bigcap_{n}[n+1,B_{n}]$. A similar argument shows that $(Q,B)$
decides $(u,b)$, for every $(u,b)\in 2^{<\infty}\times
\mathcal{P}[B]$ with $depth_{B}(b) = |u|$. Then in virtue of lemma
\ref{lem1}(f), $Q$ and $B$ are as required.

\end{proof}

\vspace{.25 cm}

\begin{lem}\label{lem3}
Let $Q$ and $B$ be as in lemma \ref{lem2}. Suppose $(Q,B)$ rejects
$(<>,\emptyset)$. Then there exists $D\leq B$ such that $(Q,D)$
rejects $(u,b)$, for every $(u,b)\in 2^{<\infty}\times
\mathcal{P}[D]$ with $depth_{D}(b)\leq |u|$.
\end{lem}
\begin{proof}

\vspace{.25 cm}

Let us build a fusion sequence $([n,D_{n}])_{n}$. Let $D_{0} = B$.
Then by hipothesis $(Q,D_{0})$ rejects $(<>,\emptyset)$. Suppose
$D_{n}$ is given such that $(Q,D_{n})$ rejects $(u,b)$, for every
$(u,b)\in 2^{n}\times \mathcal{P}[D_{n}]$ with $depth_{D_{n}}(b) =
n$.

Now, let $u_{0}$, $u_{1}$,... , $u_{2^{n+1}-1}$ be a list of the
elements of $2^{n+1}$; and let $b_{0}$, $b_{1}$,... , $b_{r}$ be a
list of the $b\in \mathcal{P}[D_{n}]$ such that $depth_{D_{n}}(b)
= n$. By lemma \ref{lem1}(f), $(Q,D_{n})$ rejects $(u_{i},b_{j})$
for every $(i,j)\in\{0, 1,\cdots ,  2^{n+1}-1\}\times \{0,
1,\cdots , r\}$. Now, by lemma \ref{lem1}(e) there exists
$D_{n}^{0,0}\in [n, D_{n}]$ such that $(Q,D_{n}^{0,0})$ rejects
$(u_{0},b)$ for every $b\in p_{|b_{0}|+1}[b_{0},D_{n}^{0,0}]$.

In the same way, for every $(i,j)\in\{0, 1,\cdots ,
2^{n+1}-1\}\times \{0, 1,\cdots , r\}$,   we can find a
$D_{n}^{i,j}$ with: $D_{n}^{i,j+1}\in [n,D_{n}^{i,j}]$,
$D_{n}^{i+1,0}\in [n,D_{n}^{i,r}]$; and such that
$(Q,D_{n}^{i,j})$ rejects $(u_{i},b)$ for every $b\in
p_{|b_{j}|+1}[b_{j},D_{n}^{i,j}]$. Let $D_{n+1} =
D_{n}^{2^{n+1}-1,r}$. Notice that if $(u,b)\in 2^{n+1}\times
\mathcal{P}[D_{n+1}]$ and $depth_{D_{n+1}}(b) = n+1$ then $u =
u_{i}$ and $b\in p_{|b_{j}|+1}[b_{j},D_{n}^{i,j}]$ for some
$(i,j)\in\{0, 1,\cdots ,  2^{n+1}-1\}\times \{0, 1,\cdots , r\}$.
Hence $(Q,D_{n+1})$ rejects $(u,b)$. Besides, $D_{n+1}\in
[n,D_{n}]$.

Now take $D\in\bigcap_{n}[n,D_{n}]$. Then $D$ is as required.

\end{proof}

\vspace{.25 cm}

\begin{proof}[Proof of theorem \ref{thm2}]

\vspace{.25 cm}

Given $\mathcal{F}\subseteq 2^{<\infty}\times\mathcal{P}$, perfect
$P\subseteq 2^{\infty}$ and $A\in\mathcal{R}$, consider the
combinatorial forcing 1. Let $Q\subseteq P$ and $B\leq A$ be as in
lemma \ref{lem2}. If $(Q,B)$ accepts $(<>,\emptyset)$ then part
(a) of theorem \ref{thm2} holds by the definition of "accepts". So
suppose $(Q,B)$ does not accept (and hence, rejects)
$(<>,\emptyset)$. By lemma \ref{lem3}, find $D\leq B$ such that
$(Q,D)$ rejects $(u,b)$, for every $(u,b)\in
2^{<\infty}\times\mathcal{P}[D]$ with $depth_{D}(b)\leq |u|$.
Suppose towards a contradiction that there exist $(t,b)$ in
$(T_{Q}\times \mathcal{P}[D])\cap\mathcal{F}$. Find $u_{t}\in
2^{<\infty}$ such that $Q(u_{t})\subseteq Q\cap [t]$. Then $(Q,D)$
accepts $(u_{t},b)$: for $x\in Q(u_{t})$ and $C\in [b,D]$, let $k
= |t|$ and $m$ be such that $p_{m}(C) = b$. Then
$(x|_{k},p_{m}(C)) = (t,b)\in \mathcal{F}$.

But then, by lemma \ref{lem1} (f), $(Q,D)$ accepts $(v,b)$, for
every $v\in 2^{<\infty}$ such that $u_{t}\sqsubseteq v$ and
$|v|\geq depth_{D}(b)$. This is a contradiction with the choice of
$D$. Therefore, for $S = Q$ and $D$ part (b) of theorem \ref{thm2}
holds.
\end{proof}

\vspace{.5 cm}

Now we are ready to prove our main result:

\vspace{.5 cm}

\begin{proof}[Proof of theorem \ref{thm1}]

\vspace{.25 cm}

\textbf{(a)} The implication from left to right is obvious. So suppose $\mathcal{X}\subseteq 2^{\infty}\times \mathcal{R}$ has the $\mathbb{P}\times Exp(\mathcal{R})$-Baire Property, and let $P\times [a,A]$ be given.  In order to make the proof notationally simpler, we will assume $a = \emptyset$ without a loss of generality.

\vspace{.25 cm}

\begin{claim} Given $\hat{\mathcal{X}}\subseteq 2^{\infty}\times
\mathcal{R}$, perfect $\hat{P}\subseteq 2^{\infty}$ and
$\hat{A}\in\mathcal{R}$, there exist $Q\subseteq \hat{P}$ and
$B\leq \hat{A}$ such that for each $(u,b)\in
2^{<\infty}\times\mathcal{P}[B]$ with $|u|\geq depth_{B}(b)$ one
of the following holds:
\begin{itemize}
\item[i.)]$Q(u)\times [b,B]\subseteq\hat{\mathcal{X}}$
\item[ii.)]$R\times [b,C]\not\subseteq\hat{\mathcal{X}}$, for
every $R\subseteq Q(u)$ and every $C\leq B$ compatible with $b$.
\end{itemize}
\end{claim}

\begin{proof}[Proof of theorem Claim]

\vspace{.25 cm}

Consider the Combinatorial Forcing 2 and apply lemma \ref{lem2}.

\end{proof}

\vspace{.25 cm}

Apply the claim to $\mathcal{X}$, $P$ and $A$ to find
$Q_{1}\subseteq P$ and $B_{1}\leq A$ such that for each $(u,b)\in
2^{<\infty}\times\mathcal{P}[B_{1}]$ with $|u|\geq
depth_{B_{1}}(b)$ one of the following holds:

\begin{itemize}
\item[(1)]$Q_{1}(u)\times [b,B_{1}]\subseteq\mathcal{X}$ or
\item[(2)]$R\times [b,C]\not\subseteq\mathcal{X}$, for every
$R\subseteq Q_{1}(u)$ and every $C\leq B_{1}$ compatible with $b$.
\end{itemize}

\vspace{.25 cm}

For each $t\in T_{Q_{1}}$, choose $u_{1}^{t}\in 2^{<\infty}$ such
that $u_{1}^{t}(Q_{1})\sqsubseteq t$. \vspace{.25 cm}

Let $$\mathcal{F}_{1} = \{(t,b)\in
T_{Q_{1}}\times\mathcal{P}[B_{1}] : Q_{1}(u_{1}^{t})\times
[b,B_{1}]\subseteq\mathcal{X}\}$$ Now, pick $S_{1}\subseteq Q_{1}$
and $D_{1}\leq B_{1}$ satisfying theorem \ref{thm2}. If (a) of
3.0.4 holds then $S_{1}\times [0,D_{1}]\subseteq \mathcal{X}$ and
we are done. So suppose (b) holds.

Apply the claim to $\mathcal{X}^{c}$, $S_{1}$ and $D_{1}$ to find
$Q_{2}\subseteq S_{1}$ and $B_{2}\leq D_{1}$ such that for each
$(u,b)\in 2^{<\infty}\times\mathcal{P}[B_{2}]$ with $|u|\geq
depth_{D_{2}}(b)$ one of the following holds:

\begin{itemize}
\item[(3)]$Q_{2}(u)\times [b,B_{2}]\subseteq\mathcal{X}^{c}$ or
\item[(4)]$R\times [b,C]\not\subseteq\mathcal{X}^{c}$, for every
$R\subseteq Q_{2}(u)$ and every $C\leq B_{2}$ compatible with $b$.
\end{itemize}

\vspace{.25 cm}

As before, for each $t\in T_{Q_{2}}$, choose $u_{2}^{t}\in
2^{<\infty}$ such that $u_{2}^{t}(Q_{2})\sqsubseteq t$.

\vspace{.25 cm}

Let $$\mathcal{F}_{2} = \{(t,b)\in
T_{Q_{2}}\times\mathcal{P}[B_{2}] : Q_{2}(u_{2}^{t})\times
[b,B_{2}]\subseteq\mathcal{X}^{c}\}$$ Again, pick $S_{2}\subseteq
Q_{2}$ and $D_{2}\leq B_{2}$ satisfying theorem \ref{thm2}. If (a)
of 3.0.4 holds then $S_{2}\times [0,D_{2}]\cap\mathcal{X} =
\emptyset$ and we are done. So suppose (b) holds again. Let us see
that this contradicts the fact that $\mathcal{X}$ has the
$\mathbb{P}\times Exp(\mathcal{R})$-Baire Property:

Notice that for every $(t,b)\in T_{S_{2}}\times\mathcal{P}[D_{2}]$ the following holds:

\vspace{.25 cm}

\begin{itemize}
\item[(i)]$Q_{1}(u_{1}^{t})\times [b,B_{1}]\not\subseteq\mathcal{X}$, and
\item[(ii)]$Q_{2}(u_{2}^{t})\times [b,B_{2}]\not\subseteq\mathcal{X}^{c}$.
\end{itemize}

\vspace{.25 cm}

So, suppose there is a nonempty $R\times [b,C]\subseteq
S_{2}\times [\emptyset,D_{2}]\cap\mathcal{X}$, and pick $t\in
T_{R}$ with $\mid u_{1}^{t}\mid \ \ \geq \ \ depth_{B_{1}}(b)$.
Notice that $R\cap [t]\subseteq Q_{1}(u_{1}^{t})$. On the one hand
we have that $R\cap [t]\times [b,C]\subseteq R\times
[b,C]\subseteq\mathcal{X}$. But in virtue of (i),
$Q_{1}(u_{1}^{t})\times [b,B_{1}]\not\subseteq\mathcal{X}$ and
hence by (2) above we have that $R\cap [t]\times
[b,C]\not\subseteq\mathcal{X}$. If we suppose that there is a
nonempty $R\times [b,C]\subseteq S_{2}\times
[\emptyset,D_{2}]\cap\mathcal{X}^{c}$ we reach to a similar
contradiction in virtue of (ii) and (4) above. So there is neither
$R\times [b,C]\subseteq S_{2}\times
[\emptyset,D_{2}]\cap\mathcal{X}$ nor $R\times [b,C]\subseteq
S_{2}\times [\emptyset,D_{2}]\cap\mathcal{X}^{c}$. But this is
impossible because $\mathcal{X}$ has the $\mathbb{P}\times
Exp(\mathcal{R})$-Baire Property.

\vspace{.25 cm}

\textbf{(b)} Again, the implication from left to right is obvious.
Conversely, the result follows easily from (a) and the fact that
$\mathcal{X}$ is $\mathbb{P}\times Exp(\mathcal{R})$-meager.

This completes the proof of theorem \ref{thm1}.

\end{proof}

\section{Some Particular Cases.}\label{cases}

Several interesting consequences can be derived  from the facts
obtained in the previous section. Some of them are known classical
results and the others are parametrized versions of known
Ellentuck type theorems in nonclassical spaces.

\vspace{1.25 cm}

Let $k$ be a positive integer. For $p:\mathbb{N}\rightarrow\{0,1, \cdots, k\}$, let $supp(p)$ denote the set $\{n : p(n)\neq 0\}$ and let $rang(p)$ denote the range of $p$.

\vspace{.25 cm}

Let us consider the set $$FIN_{k} := \{p:\mathbb{N}\rightarrow\{0,1, \cdots, k\} : supp(p)\ \ \mbox{is finite and}\ \ k\in rang(p)\}.$$

A {\em block basic sequence} is any finite or infinite sequence $X = (x_{n})_{n\in I\subseteq\mathbb{N}}$ of elements of $FIN_{k}$ such that $$max(supp(x_{n})) < min(supp(x_{m}))\ \ \mbox{whenever}\ \ n < m.$$ We shall use $a$, $b$, $c$, ... for finite block basic sequences, and $A$, $B$, $C$, ... for infinite block basic sequences. In this latter case we will assume that the set of indexes is $I = \mathbb{N}$.

\vspace{.25 cm}

Define $T : FIN_{k}\rightarrow FIN_{k-1}$ by $$T(p)(n)=
max\{p(n)-1, 0\}.$$ In \cite{Tod} this is called the {\em tetris
operation}. For every $j\in\mathbb{N}$, $T^{(j)}$ is the j-th
iteration of $T$, where $T^{(0)}(p) = p$.

\vspace{.25 cm}

For a given block basic sequence $X = (x_{n})_{n\in I\subseteq\mathbb{N}}$, the {\em subspace} of $FIN_{k}$ {\em generated} by $X$, denoted by $[X]$, is the set of elements of $FIN_{k}$ of the form: $$T^{(j_{0})}(x_{n_{0}})+T^{(j_{1})}(x_{n_{1}})+\cdots+T^{(j_{r})}(x_{n_{r}})$$ where $n_{0} < n_{1} < \cdots < n_{r}$ is a finite sequence of elements of $I$ and $j_{0} < j_{1} < \cdots < j_{r}$ is a sequence of elements of $\{0,1, \cdots, k\}$ such that $j_{i} = 0$ for some $i\in \{0,1, \cdots, r\}$.

The following result shows us an important feature of $FIN_{k}$, which provides us of a pigeon hole principle within this context:
\vspace{.25 cm}

\begin{thm}\label{thmGow}
[Gowers \cite{Gow}] For every integer $r > 0$ and every partition
$\phi : FIN_{k}\rightarrow\{0, 1, \cdots, r-1\}$ there exists an
infinite block basic sequence $A$ such that $\phi$ is constant in
$[A]$.
\end{thm}
\qed

\vspace{.25 cm}

In the case $k = 1$, $FIN_{k}$ is the set $FIN$ of nonempty finite
subsets of $\mathbb{N}$, and theorem \ref{thmGow} is Hindman's
theorem \cite{Hin}.

\vspace{.25 cm}

Let $FIN_{k}^{[\infty]}$ be the set of infinite block basic sequences and define,  $$A\leq B\ \ \mbox{iff}\ \ A\subseteq [B]$$ for $A,B\in FIN_{k}^{[\infty]}$. Also, for every $A\in FIN_{k}^{[\infty]}$, let  the n-th approximation of $A$ be $$p_{n}(A) = \ \ \mbox{the first}\ \ n\ \ \mbox{elements of}\ \ A.$$

\vspace{.25 cm}

Then the set $\mathcal{P}$ of approximations is
$FIN_{k}^{[<\infty]}$, the set of finite block basic sequences.
Now, for $a,b\in FIN_{k}^{[<\infty]}$ define $a\leq_{fin} b$ if
and only if $$a=b=\emptyset\ \ \mbox{or}\ \ a\subset [b]\ \
\mbox{and}\ \ max(supp\bigcup a) = max(supp\bigcup b).$$ With this
terminology and the obvious definition of the neighborhoods
$[a,A]$ (and the family $Exp(FIN_{k}^{[\infty]})$), the triplet
$(FIN_{k}^{[\infty]}, \leq, (p_{n})_{n\in\mathbb{N}})$ satisfies
A1-A6. Also, $FIN_{k}^{[\infty]}$ is readily identified to a
closed subset of $\mathcal{P}^{[\infty]}$ and hence it is a Ramsey
space, in virtue of the Abstract Ellentuck Theorem (\ref{AbsEll}
above). Here A6 reduces to a natural variation of Gowers' theorem
(\ref{thmGow} above). This yields the following corollary of
theorem \ref{thm1}:

\begin{cor}\label{ParMil}
Let $\mathcal{X}\subseteq 2^{\infty}\times FIN_{k}^{[\infty]}$ be
given. $\mathcal{X}$ is perfectly Ramsey\ \ iff\ \ $\mathcal{X}$
has the $\mathbb{P}\times Exp(FIN_{k}^{[\infty]})$-Baire Property.
$\mathcal{X}$ is perfectly Ramsey null\ \ iff\ \ $\mathcal{X}$ is
$\mathbb{P}\times Exp(FIN_{k}^{[\infty]})$-meager.
\end{cor}

\qed

\vspace{.25 cm}

In the same way, in virtue of corollary \ref{Ell}, we have the
following:

\vspace{.25 cm}

\begin{cor}\label{ParEll}
Let $\mathcal{X}\subseteq 2^{\infty}\times\mathbb{N}^{[\infty]}$ be given. $\mathcal{X}$ is perfectly Ramsey\ \ iff\ \ $\mathcal{X}$ has the $\mathbb{P}\times Exp(\mathbb{N}^{[\infty]})$-Baire Property. $\mathcal{X}$ is perfectly Ramsey null\ \ iff\ \ $\mathcal{X}$ is $\mathbb{P}\times Exp(\mathbb{N}^{[\infty]})$-meager.
\end{cor}

\qed

\vspace{.25 cm}

Corollary \ref{ParEll} is the parametrization of Ellentuck's
theorem \cite{Ell} obtained by Pawlikowski in \cite{Paw}. And
corollary \ref{ParMil} gives a parametrized version of Milliken's
theorem \cite{Milk}, when  $k = 1$, and a parametrized version of
the corresponding Ellentuck type theorem in the Ramsey space
$(FIN_{k}^{[\infty]}, \leq, (p_{n})_{n\in\mathbb{N}})$ defined in
this section, when $k
> 1$.

\vspace{.25 cm}

\section{Closedness Under the Souslin Operation.}

In this section we go back to our parametrization of (abstract) Ramsey spaces and study the family of perfectly Ramsey sets in relation to the Souslin operation. First, we proof the following fact which turns out to be crucial for this study.

\vspace{.25 cm}

\begin{lem}\label{RamseyNull}
Given $(\mathcal{R}, \leq, (p_{n})_{n\in\mathbb{N}})$ satisfying
A1-A6 and with $\mathcal{R}$ metrically closed, the perfectly
Ramsey null subsets of $2^{\infty}\times\mathcal{R}$ form a
$\sigma$-ideal.
\end{lem}
\begin{proof}

\vspace{.25 cm}

Let $(\mathcal{X}_n)_n$ be a sequence of perfectly Ramsey null
subsets of $2^{\infty}\times\mathcal{R}$ and fix $P\times [a,A]$.
We can assume $a = \emptyset$. Also notice that the finite union
of perfectly Ramsey null sets yields a perfectly Ramsey null set;
so we will assume $(\forall n)\ \mathcal{X}_n\subseteq
\mathcal{X}_{n+1}$. Proceeding as in the proof of lemma \ref{lem2}
we build fusion sequences $Q_n$, $[n+1 ,B_n]$ as follows: take
$Q_0\subseteq P$, $B_0\leq A$ such that $Q_0\times [0,B_0]\cap
\mathcal{X}_0 = \emptyset$. Suppose $Q_n$, $[n+1 ,B_n]$ have been
defined such that
$$Q_n\times [b,B_n]\cap \mathcal{X}_n = \emptyset$$

\vspace{.25 cm} for every $b\in\mathcal{P}[B_n]$ with
$depth_{B_n}(b) = n$. Since $\mathcal{X}_{n+1}$ is perfectly
Ramsey null, applying this fact successively we find
$Q_{n+1}\subseteq_{n+1} Q_n$, and $B_{n+1}\in [n+1 ,B_n]$ such
that
$$Q_{n+1}\times [b,B_{n+1}]\cap \mathcal{X}_{n+1} = \emptyset$$ for every $b\in\mathcal{P}[B_{n+1}]$ with $depth_{B_{n+1}}(b) =
n+1$. Let $$Q = \bigcap_n Q_n\ \ \mbox{and}\ \ B = \bigcap_n [n+1
,B_n].$$ Then $Q\times [0,B]\cap \bigcup_n\mathcal{X}_n =
\emptyset$: take $(x,C)\in Q\times [0,B]$ and fix arbitrary $n$.
To show that $(x,C)\not\in\mathcal{X}_n$ let $k$ be large enough
so that $depth_B(p_k(C)) = m \geq n$. Then by construction
$Q\times [p_k(C),B]\cap \mathcal{X}_m = \emptyset$ and hence,
since $\mathcal{X}_n\subseteq \mathcal{X}_m$, we have
$(x,C)\not\in\mathcal{X}_n$. This completes the proof.

\vspace{.25 cm}

\end{proof}

Now, we borrow some terminology from \cite{Paw}:

\vspace{.25 cm}

Let $\mathcal{A}$ be a family of subsets of a set $\mathcal{Z}$.
We say that $\mathcal{X}, \mathcal{Y}\subseteq\mathcal{Z}$ are
{\em compatible} (with respect to $\mathcal{A}$) if there exists
$\mathcal{W}\in \mathcal{A}$ such that
$\mathcal{W}\subseteq\mathcal{X}\cap\mathcal{Y}$. Also, we say
that $\mathcal{A}$ is $M$-{\em like} if for any
$\mathcal{B}\subseteq\mathcal{A}$  such that $|\mathcal{B}| <
|\mathcal{A}|$, every member of  $\mathcal{A}$ which is not
compatible with any member of  $\mathcal{B}$ is compatible with
$\mathcal{Z}\setminus\bigcup\mathcal{B}$.

\vspace{1 cm}

Notice that the family $\mathbb{P}$ of perfect subsets of
$2^{\infty}$ is $M$-like, as well as the family $Exp(\mathcal{R})$
(this is true of any topological basis). Therefore, according to
lemma 2.7 in \cite{Paw}, if we require that $|Exp(\mathcal{R})| =
|\mathbb{P}|$ (= $2^{\aleph_0}$), then the family
$\mathbb{P}\times Exp(\mathcal{R}) = \{P\times [n,A] :
P\in\mathbb{P}\ \mbox{and}\ \   A\in\mathcal{R}\}$ is also
$M$-like. This lead us to the following:

\vspace{.25 cm}

\begin{cor}
Let $(\mathcal{R}, \leq, (p_{n})_{n\in\mathbb{N}})$ satisfying
A1-A6, with $\mathcal{R}$ metrically closed be such that
$|\mathcal{R}| = 2^{\aleph_0}$. Then, the family of perfectly
Ramsey subsets of $2^{\infty}\times\mathcal{R}$ is closed under
the Souslin operation.
\end{cor}

\begin{proof}

\vspace{.25 cm}

In virtue of theorem \ref{thm1}, the family of perfectly Ramsey
subsets of $2^{\infty}\times\mathcal{R}$ coincides with the family
of subsets of $2^{\infty}\times\mathcal{R}$ which have the
$\mathbb{P}\times Exp(\mathcal{R})$-Baire property. And as we
pointed out in the previous parragraph, $\mathbb{P}\times
Exp(\mathcal{R})$ is $M$-like. So the proof follows from lemma
\ref{RamseyNull}on top of this section and lemmas 2.5 and 2.6 of
\cite{Paw} (which refer to a well-known result of Marczewski
\cite{Mar}).

\vspace{.25 cm}

\end{proof}

\vspace{.25 cm}

\begin{cor}
(a)[Pawlikowski] The family of perfectly Ramsey subsets of
$2^{\infty}\times\mathbb{N}^{[\infty]}$ is closed under the
Souslin operation. (b)The family of perfectly Ramsey subsets of
$2^{\infty}\times FIN_{k}^{[\infty]}$ is closed under the Souslin
operation.
\end{cor}

\qed

\begin{cor}
(a)[Miller-Todorcevic] Analytic subsets of
$2^{\infty}\times\mathbb{N}^{[\infty]}$ are perfectly Ramsey.
(b)Analytic subsets of $2^{\infty}\times FIN_{k}^{[\infty]}$ are
perfectly Ramsey.
\end{cor}

\qed

\vspace{2 cm}

{\bf Acknowledgement.} The author would like to thank Carlos Di
Prisco for many years of guidance and teachings, Stevo Todorcevic
for insightful sugerences, and Mar\' ia Carrasco and Franklin
Galindo for unvaluable feedback.

\end{document}